\documentclass[11pt]{article}
\usepackage{amssymb}

\usepackage{amsmath}
\usepackage{theorem}
\usepackage{here}
\usepackage[dvipdfmx]{color}

\newtheorem{prp}{Proposition}[section]
\theorembodyfont{\rmfamily}
\newtheorem{algo}{Algorithm}[section]
\makeatletter

\@addtoreset{equation}{section}
\makeatother

\usepackage{comment} %
\usepackage{bm}
\usepackage[dvipdfmx]{graphicx} %

\def\ep{{\varepsilon}}

\def\bbe{{\text{\boldmath $\beta$}}}

\def\bep{{\text{\boldmath $\varepsilon$}}}

\def\bth{{\text{\boldmath $\theta$}}}

\def\bmu{{\text{\boldmath $\mu$}}}

\def\Th{{\Theta}}

\def\Si{{\Sigma}}
\def\Ga{{\Gamma}}
\def\Om{{\Omega}}

\def\bTh{{\text{\boldmath $\Th$}}}

\def\bSi{{\text{\boldmath $\Si$}}}
\def\bGa{{\text{\boldmath $\Ga$}}}
\def\bOm{{\text{\boldmath $\Om$}}}

\def\bXi{{\text{\boldmath $\Xi$}}}

\def\bPsi{{\text{\boldmath $\Psi$}}}

\def\u{{\text{\boldmath $u$}}}

\def\w{{\text{\boldmath $w$}}}
\def\x{{\text{\boldmath $x$}}}
\def\y{{\text{\boldmath $y$}}}

\def\A{{\text{\boldmath $A$}}}
\def\B{{\text{\boldmath $B$}}}
\def\C{{\text{\boldmath $C$}}}

\def\I{{\text{\boldmath $I$}}}

\def\M{{\text{\boldmath $M$}}}

\def\O{{\text{\boldmath $O$}}}
\def\P{{\text{\boldmath $P$}}}
\def\Q{{\text{\boldmath $Q$}}}

\def\U{{\text{\boldmath $U$}}}

\def\X{{\text{\boldmath $X$}}}
\def\Y{{\text{\boldmath $Y$}}}
\def\Z{{\text{\boldmath $Z$}}}

\def\Xbt{{\widetilde \X}}
\def\Ybt{{\widetilde \Y}}

\def\tr{{\rm tr\,}}

\def\etr{{\rm etr\,}}

\def\[{{\text{\boldmath $[$}}}
\def\]{{\text{\boldmath $]$}}}

\def\/{{\Bigr/\!\!}}

\def\1r{{\rm (1)}}
\def\2r{{\rm (2)}}
\def\3r{{\rm (3)}}
\def\4r{{\rm (4)}}
\def\5r{{\rm (5)}}

\def\non{{\nonumber}}
\DeclareMathOperator*{\bdiag}{{\bf Diag\,}}
\def\bBe{{\text{\boldmath $\mathcal{B}$}}}

\begin{document}
\title{A Short Note on the Efficiency of Markov Chains for Bayesian Linear Regression Models with Heavy-Tailed Errors}
\author{
Yasuyuki Hamura\footnote{Graduate School of Economics, Kyoto University, 
Yoshida-Honmachi, Sakyo-ku, Kyoto, 606-8501, JAPAN. 
\newline{
E-Mail: yasu.stat@gmail.com}} \
}
\maketitle
\begin{abstract}
In this short note, we consider posterior simulation for a linear regression model when the error distribution is given by a scale mixture of multivariate normals. 
We first show that the sampler of Backlund and Hobert (2020) for the case of the conditionally conjugate normal-inverse Wishart prior continues to be geometrically ergodic even when the error density is heavier-tailed. 
Moreover, we prove that the ergodicity is uniform by verifying the minorization condition. 
In the second half of this note, we treat an improper case and show that the sampler of Section 4 of Roy and Hobert (2010) is geometrically ergodic under significantly milder conditions.

\par\vspace{4mm}
\end{abstract}

\section{Introduction}
\label{sec:introduction}
Scale mixtures of normals, including Student's $t$ distribution, are widely used as convenient error densities for the traditional regression modeling of data contaminated with outliers. 
In the Bayesian framework, it is difficult to directly sample from the marginal posterior distribution under a regression model with a mixture error density but we can introduce latent or mixing variables to construct simple Gibbs algorithms to obtain posterior samples of location and scale parameters; see, for example, Hobert et al. (2018) and the references therein. 
Although such a data augmentation approach is taken by many researchers, efficiency properties of simulation algorighms have not been fully investigated from a theoretical point of view and, in recent years also, papers have been published in this area; see, for example, the references cited in Ekvall and Jones (2021). 

A popular efficiency criterion is geometric ergodicity. 
It is not easy to establish in general but we can sometimes use the drift and minorization technique to prove it; see Jones and Hobert (2001). 
Much work has been done to show geometric ergodicity for regression models with mixture error densities. 
For example, Roy and Hobert (2010) introduced a useful energy function in a simple case. 
Hobert et al. (2018) considered the case of relatively light-tailed error density and improper prior distribution. 
Jung and Hobert (2014) and Qin and Hobert (2018) derived a property stronger than geometric ergodicity under improper priors. 
Backlund and Hobert (2020) extended the result of Hobert et al. (2018) to the of conditionally conjugate proper priors. 
Backlund et al. (2021) considered the case of independent proper priors. 
However, the case of a heavier-tailed error density and a conditionally conjugate normal-inverse Wishart prior has not been studied, to the author's %
limited knowledge of this area. 

In the first half of this short note, we consider posterior simulation for a linear regression model with a conditionally conjugate normal-inverse Wishart prior when the error distribution is given by a scale mixture of multivariate normals. 
We show that the sampler of Backlund and Hobert (2020) continues to be geometrically ergodic even when the error density is heavier-tailed. 
Moreover, we prove that it is actually uniformly ergodic by verifying the minorization condition. 
Next, we turn to the case of an improper prior. 
Using a different energy function, we show that the sampler of Section 4 of Roy and Hobert (2010) is geometrically ergodic under significantly milder conditions.

\section{The Model and Algorithm}
\label{sec:Po} 
Consider the linear regression model of Backlund and Hobert (2020). 
Specifically, let $h \colon (0, \infty ) \to [0, \infty )$ be a normalized mixing density and suppose that 
\begin{align}
&\y _i \sim {1 \over | \bSi |^{1 / 2}} f_h ( \bSi ^{- 1 / 2} ( \y _i - \bBe ^{\top } \x _i )) \non 
\end{align}
for $i = 1, \dots , n$, where $\Y = ( \y _1 , \dots , \y _n )^{\top } \in \mathbb{R} ^{n \times d}$ and $\X = ( \x _1 , \dots , \x _n )^{\top } \in \mathbb{R} ^{n \times p}$ are outcome and explanatory variables, $\bBe \in \mathbb{R} ^{p \times d}$ and $\bSi > \O ^{(d)}$ are the matrix of regression coefficients and the covariance matrix, and where 
\begin{align}
f_h ( \bep ) &= \int_{0}^{\infty } {u^{d / 2} \over (2 \pi )^{d / 2}} \exp \Big( - {u \over 2} \| \bep \| ^2 \Big) h(u) du \text{,} \quad \bep \in \mathbb{R} ^d \text{,} \non 
\end{align}
is the error density. 
Generally, mixing densities with heavier left tails induce error densities with heavier tails. 
The prior distribution is given by 
\begin{align}
p( \bBe | \bSi ) p( \bSi ) &= {\rm{N}}_{p, d} ( \bBe | \B , \A , \bSi ) {\rm{IW}}_d ( \bSi | \nu , \bTh ) \text{,} \non 
\end{align}
where $\B \in \mathbb{R} ^{p \times d}$, $\A > \O ^{(p)}$, $\nu > d - 1$, and $\bTh > \O ^{(d)}$ are hyperparameters, and it is proper and conditionally conjugate. 
We use the notation $\B $ instead of $\bth $ to clarify that it is a matrix. 

As shown in Backlund and Hobert (2020), although the marginal posterior density $p( \bBe , \bSi | \Y )$ is intractable, we can introduce augmenting variables $\u = ( u_i )_{i = 1}^{n} \in (0, \infty )^n$ to construct the following Gibbs algorithm to sample $( \bBe , \bSi )$. 

\begin{algo}[Backlund and Hobert (2020)]
\label{algo:DA} 
The parameters $\bBe $ and $\bSi $ are updated in the following way. 
\begin{itemize}
\item
For each $i = 1, \dots , n$, sample $u_i \sim p( u_i | \bBe , \bSi , \Y )$, where 
\begin{align}
p( u_i | \bBe , \bSi , \y ) &\propto h( u_i ) {u_i}^{d / 2} \exp \{ - ( u_i / 2) ( \y _i - \bBe ^{\top } \x _i )^{\top } \bSi ^{- 1} ( \y _i - \bBe ^{\top } \x _i ) \} \text{.} \non 
\end{align}
Let $\U = \bdiag \u > \O ^{(n)}$. 
\item
Sample $\bSi \sim {\rm{IW}}_d (n + \nu , \bPsi ^{- 1} )$ and then $\bBe \sim {\rm{N}}_{p, d} ( \bGa , \bOm , \bSi )$, where 
\begin{align}
&\bPsi = \bTh ^{- 1} + \B ^{\top } \A ^{- 1} \B + \Y ^{\top } \U \Y - \bGa ^{\top } \bOm ^{- 1} \bGa \text{,} \non \\
&\bGa = ( \X ^{\top } \U \X + \A ^{- 1} )^{- 1} ( \X ^{\top } \U \Y + \A ^{- 1} \B ) \text{,} \quad \text{and} \non \\
&\bOm = ( \X ^{\top } \U \X + \A ^{- 1} )^{- 1} \text{.} \non 
\end{align}
\end{itemize}
\end{algo}

\section{Geometric Ergodicity}
\label{sec:result} 
The following result shows that Algorithm \ref{algo:DA} is efficient. 

\begin{prp}
\label{prp:ergodicity} 
Suppose that $\int_{0}^{\infty } u^{d / 2 + 2} h(u) du < \infty $. 
Then the Markov chain based on Algorithm \ref{algo:DA} is geometrically ergodic. 
\end{prp}

This result is different from Proposition 1 of Backlund and Hobert (2020) in that we can apply Proposition \ref{prp:ergodicity} when the mixing density $h$ is heavy-tailed near the origin. 
An important special case is that of a gamma mixing density with a small shape parameter, which corresponds to Student's $t$ error density with small degrees of freedom. 

\bigskip

\noindent
{\bf Proof of Proposition \ref{prp:ergodicity}.} \ \ Many steps of our proof is same as those of the proof of Proposition 1 of Backlund and Hobert (2020). 
The only additional key idea is to exploit the fact that when the prior is proper, %
there are positive constants in the joint density, such as rate hyperparameters, which could be useful in bounding some terms. 
In this proof, we suppress the dependence on $\Y $. 
Let 
\begin{align}
V( \bBe , \bSi ) &= \sum_{i = 1}^{n} ( \y _i - \bBe ^{\top } \x _i )^{\top } \bSi ^{- 1} ( \y _i - \bBe ^{\top } \x _i ) \text{.} \non 
\end{align}
As in Hobert et al. (2018), the minorization condition associated with this energy function holds since $\int_{0}^{\infty } u^{d / 2} h(u) du < \infty $ by assumption. 

By calculations on pages 5 and 6 of Backlund and Hobert (2020), 
\begin{align}
&E[ V( \bBe , \bSi ) | \u ] = \sum_{i = 1}^{n} ( \y _i - \bGa ^{\top } \x _i )^{\top } E[ \bSi ^{- 1} | \u ] ( \y _i - \bGa ^{\top } \x _i ) + d \sum_{i = 1}^{n} {\x _i}^{\top } \bOm \x _i \text{.} \non 
\end{align}
where $E[ \bSi ^{- 1} | \u ] = (n + \nu ) ( \bTh ^{- 1} + \B ^{\top } \A ^{- 1} \B + \Y ^{\top } \U \Y - \bGa ^{\top } \bOm ^{- 1} \bGa )^{- 1}$. 
The posterior is proper regardless of the value of $\bTh ^{- 1}$ and we have $\B ^{\top } \A ^{- 1} \B + \Y ^{\top } \U \Y - \bGa ^{\top } \bOm ^{- 1} \bGa \ge \O ^{(d)}$. 
(Indeed, 
\begin{align}
\bGa ^{\top } \bOm ^{- 1} \bGa &= ( \X ^{\top } \U \Y + \A ^{- 1} \B )^{\top } \{ \A - \A \X ^{\top } ( \U ^{- 1} + \X \A \X ^{\top } )^{- 1} \X \A \} ( \X ^{\top } \U \Y + \A ^{- 1} \B ) \text{,} \non 
\end{align}
where 
\begin{align}
&\B ^{\top } \A ^{- 1} \{ \A - \A \X ^{\top } ( \U ^{- 1} + \X \A \X ^{\top } )^{- 1} \X \A \} \X ^{\top } \U \Y \non \\
&= \B ^{\top } \X ^{\top } \U \Y - \B ^{\top } \X ^{\top } ( \U ^{- 1} + \X \A \X ^{\top } )^{- 1} (- \U ^{- 1} + \U ^{- 1} + \X \A \X ^{\top } ) \U \Y \non \\
&= \B ^{\top } \X ^{\top } ( \U ^{- 1} + \X \A \X ^{\top } )^{- 1} \Y \non 
\end{align}
and similarly 
\begin{align}
&\Y ^{\top } \U \X \{ \A - \A \X ^{\top } ( \U ^{- 1} + \X \A \X ^{\top } )^{- 1} \X \A \} \X ^{\top } \U \Y = \Y ^{\top } \U \Y - \Y ^{\top } ( \U ^{- 1} + \X \A \X ^{\top } )^{- 1} \Y \text{,} \non 
\end{align}
and it follows that 
\begin{align}
&\B ^{\top } \A ^{- 1} \B + \Y ^{\top } \U \Y - \bGa ^{\top } \bOm ^{- 1} \bGa %
= ( \X \B - \Y )^{\top } ( \U ^{- 1} + \X \A \X ^{\top } )^{- 1} ( \X \B - \Y ) \text{,} \non 
\end{align}
which is nonnegative definite.) 
Thus, 
\begin{align}
E[ V( \bBe , \bSi ) | \u ] &\le (n + \nu ) \sum_{i = 1}^{n} ( \y _i - \bGa ^{\top } \x _i )^{\top } \bTh ( \y _i - \bGa ^{\top } \x _i ) + d \sum_{i = 1}^{n} {\x _i}^{\top } \bOm \x _i \non \\
&\le 2 (n + \nu ) \sum_{i = 1}^{n} ( {\y _i}^{\top } \bTh \y _i + {\x _i}^{\top } \bGa \bTh \bGa ^{\top } \x _i ) + d \sum_{i = 1}^{n} {\x _i}^{\top } \A \x _i \non \\
&= d \sum_{i = 1}^{n} {\x _i}^{\top } \A \x _i + 2 (n + \nu ) \sum_{i = 1}^{n} {\y _i}^{\top } \bTh \y _i + 2 (n + \nu ) \sum_{i = 1}^{n} \tr ( {\x _i}^{\top } \bGa \bTh \bGa ^{\top } \x _i ) \text{,} \non 
\end{align}
where the first two terms on the right side are constants. 

Now, let $M_1 > 0$ be such that $\X ^{\top } \X \le M_1 \A ^{- 1}$. 
Then 
\begin{align}
\sum_{i = 1}^{n} \tr ( {\x _i}^{\top } \bGa \bTh \bGa ^{\top } \x _i ) &= \tr \Big\{ \sum_{i = 1}^{n} ( \bTh ^{1 / 2} )^{\top } \bGa ^{\top } \x _i {\x _i}^{\top } \bGa \bTh ^{1 / 2} \Big\} \le M_1 \tr \{ ( \bTh ^{1 / 2} )^{\top } \bGa ^{\top } \A ^{- 1} \bGa \bTh ^{1 / 2} \} \text{.} \non 
\end{align}
We have 
\begin{align}
\bGa ^{\top } \A ^{- 1} \bGa &= ( \X ^{\top } \U \Y + \A ^{- 1} \B )^{\top } ( \A ^{- 1} + 2 \X ^{\top } \U \X + \X ^{\top } \U \X \A \X ^{\top } \U \X )^{- 1} ( \X ^{\top } \U \Y + \A ^{- 1} \B ) \non \\
&\le ( \X ^{\top } \U \Y + \A ^{- 1} \B )^{\top } \A ( \X ^{\top } \U \Y + \A ^{- 1} \B ) \le 2 \Y ^{\top } \U ^{\top } \X \A \X ^{\top } \U \Y + 2 \B ^{\top } \A ^{- 1} \B \non \\
&\le 2 \B ^{\top } \A ^{- 1} \B + 2 M_2 \Y ^{\top } \U ^2 \Y \le 2 \B ^{\top } \A ^{- 1} \B + 2 M_2 \tr ( \U ^2 ) \Y ^{\top } \Y \non 
\end{align}
for some $M_2 > 0$. 
(These inequalities are used also in Hamura (2024).) 
Therefore, 
\begin{align}
E[ V( \bBe , \bSi ) | \u ] &\le M_3 + M_4 \tr ( \U ^2 ) \non 
\end{align}
for some $M_3 , M_4 > 0$. 
Since 
\begin{align}
E[ {u_i}^2 | \bBe , \bSi ] &= \frac{ \displaystyle \int_{0}^{\infty } h( u_i ) {u_i}^{d / 2} {u_i}^2 \exp \{ - ( u_i / 2) ( \y _i - \bBe ^{\top } \x _i )^{\top } \bSi ^{- 1} ( \y _i - \bBe ^{\top } \x _i ) \} d{u_i} }{ \displaystyle \int_{0}^{\infty } h( u_i ) {u_i}^{d / 2} \exp \{ - ( u_i / 2) ( \y _i - \bBe ^{\top } \x _i )^{\top } \bSi ^{- 1} ( \y _i - \bBe ^{\top } \x _i ) \} d{u_i} } \non \\
&\le \int_{0}^{\infty } h( u_i ) {u_i}^{d / 2} {u_i}^2 d{u_i} / \int_{0}^{\infty } h( u_i ) {u_i}^{d / 2} d{u_i} \non 
\end{align}
for all $i = 1, \dots , n$ by the covariance inequality, it follows that 
\begin{align}
E[ V( \bBe , \bSi ) | \u ] &\le M_5 \non 
\end{align}
for some $M_5 > 0$. 
The drift condition holds and this completes the proof. 
\hfill$\Box$

\section{Uniform Ergodicity}
\label{sec:unif} 
If $n \ge p$ and $\X $ is of full rank, a stronger conclusion holds.

\begin{prp}
\label{prp:unif} 
Assume that $\X ^{\top } \X $ is of full rank. 
Suppose that $\int_{0}^{\infty } u^{d / 2} h(u) du < \infty $. 
Then the Markov chain based on Algorithm \ref{algo:DA} is uniformly ergodic. 
\end{prp}

\bigskip

\noindent
{\bf Proof of Proposition \ref{prp:unif}.} \ \ In order to apply Theorem 1 of Jones and Hobert (2001), we show that there exist a proper density $q$ for $( \bBe , \bSi )$ and a positive constant $\ep > 0$ such that for all $\bBe , \bBe ^{\rm{o}} \in \mathbb{R} ^{p \times d}$ and all $\bSi , \bSi ^{\rm{o}} > \O ^{(d)}$, we have 
\begin{align}
&E^{\u | ( \bBe , \bSi )} [ p( \bBe , \bSi | \u ) | \bBe ^{\rm{o}} , \bSi ^{\rm{o}} ] \ge \ep q( \bBe , \bSi ) \text{.} \non 
\end{align}
By the proof of Proposition 1 of Backlund and Hobert (2020), $p( \bBe | \bSi , \u ) = {\rm{N}}_{p, d} ( \bBe | \bGa , \bOm , \bSi )$ and $p( \bSi | \u ) = {\rm{IW}}_d ( \bSi | n + \nu , \bPsi ^{- 1} )$. 
By Appendix A of Backlund and Hobert (2020), 
\begin{align}
p( \bbe , \bSi | \u ) &= c_0 {| \bPsi |^{(n + \nu ) / 2} / | \bOm |^{d / 2} \over | \bSi |^{(n + \nu + d + 1 + p) / 2}} \etr (- \bPsi \bSi ^{- 1} / 2) \etr \{ - \bOm ^{- 1} ( \bBe - \bGa ) \bSi ^{- 1} ( \bBe - \bGa )^{\top } / 2 \} \text{,} \non 
\end{align}
where $1 / c_0 = (2 \pi )^{p d / 2} 2^{(n + \nu ) d / 2} \pi ^{d (d - 1) / 4} \prod_{g = 1}^{d} \Ga (( n + \nu + 1 - g) / 2)$. 
Note that $| \bPsi | \ge | \bTh ^{- 1} |$ by the proof of Proposition \ref{prp:ergodicity} and that $| \bOm | \le | \A |$ as in the proof of Proposition 3.1 of Choi and Hobert (2013). 
Then 
\begin{align}
p( \bbe , \bSi | \u ) &\ge \{ c_1 / | \bSi |^{(n + \nu + d + 1 + p) / 2} \} \etr (- \bPsi \bSi ^{- 1} / 2) \etr \{ - \bOm ^{- 1} ( \bBe - \bGa ) \bSi ^{- 1} ( \bBe - \bGa )^{\top } / 2 \} \text{,} \non 
\end{align}
where $c_1 = c_0 | \bTh ^{- 1} |^{(n + \nu ) / 2} / | \A |^{d / 2}$. 
Note that 
\begin{align}
\tr ( \bPsi \bSi ^{- 1} ) &= \tr \{ \bSi ^{- 1 / 2} \bPsi ( \bSi ^{- 1 / 2} )^{\top } \} \le \tr \{ \bSi ^{- 1 / 2} ( \bTh ^{- 1} + \B ^{\top } \A ^{- 1} \B + \Y ^{\top } \U \Y ) ( \bSi ^{- 1 / 2} )^{\top } \} \text{.} \non 
\end{align}
Then 
\begin{align}
p( \bbe , \bSi | \u ) &\ge \{ c_1 / | \bSi |^{(n + \nu + d + 1 + p) / 2} \} \etr \{ - \bSi ^{- 1} ( \bTh ^{- 1} + \B ^{\top } \A ^{- 1} \B + \Y ^{\top } \U \Y ) / 2 \} \non \\
&\quad \times \etr \{ - \bSi ^{- 1} ( \bBe - \bGa )^{\top } \bOm ^{- 1} ( \bBe - \bGa ) / 2 \} \non \\
&\ge \{ c_1 / | \bSi |^{(n + \nu + d + 1 + p) / 2} \} \etr \{ - \bSi ^{- 1} ( \bTh ^{- 1} + \B ^{\top } \A ^{- 1} \B + \Y ^{\top } \U \Y ) / 2 \} \non \\
&\quad \times \etr \{ - (2 \bSi ^{- 1} \bBe ^{\top } \bOm ^{- 1} \bBe + 2 \bSi ^{- 1} \bGa ^{\top } \bOm ^{- 1} \bGa ) / 2 \} \non \\
&= \{ c_1 / | \bSi |^{(n + \nu + d + 1 + p) / 2} \} \etr \{ - \bSi ^{- 1} ( \bTh ^{- 1} + \B ^{\top } \A ^{- 1} \B + 2 \bBe ^{\top } \A ^{- 1} \bBe ) / 2 \} \non \\
&\quad \times \etr \{ - \bSi ^{- 1} (2 \bBe ^{\top } \X ^{\top } \U \X \bBe + \Y ^{\top } \U \Y ) / 2 \} \etr (- 2 \bSi ^{- 1} \bGa ^{\top } \bOm ^{- 1} \bGa / 2) \text{,} \non 
\end{align}
where the second inequality follows since 
\begin{align}
\y ^{\top } ( \bBe - \bGa )^{\top } \bOm ^{- 1} ( \bBe - \bGa ) \y &\le 2 \y ^{\top } ( \bBe ^{\top } \bOm ^{- 1} \bBe + \bGa ^{\top } \bOm ^{- 1} \bGa ) \y \non 
\end{align}
for all $\y \in \mathbb{R} ^d$. 

Now, 
\begin{align}
\bGa ^{\top } \bOm ^{- 1} \bGa &= ( \X ^{\top } \U \Y + \A ^{- 1} \B )^{\top } ( \X ^{\top } \U \X + \A ^{- 1} )^{- 1} ( \X ^{\top } \U \Y + \A ^{- 1} \B ) \non \\
&\le 2 ( \A ^{- 1} \B )^{\top } ( \X ^{\top } \U \X + \A ^{- 1} )^{- 1} \A ^{- 1} \B + 2 ( \X ^{\top } \U \Y )^{\top } ( \X ^{\top } \U \X + \A ^{- 1} )^{- 1} \X ^{\top } \U \Y \non \\
&\le 2 \B ^{\top } \A ^{- 1} \B + 2 \Y ^{\top } \U \X ( \X ^{\top } \U \X )^{- 1} ( \U \X )^{\top } \Y \text{,} \non 
\end{align}
where $\X ^{\top } \U \X \ge ( \min_{1 \le i \le n} u_i ) \X ^{\top } \X > \O ^{(p)}$ by assumption, and 
\begin{align}
&\U \X ( \X ^{\top } \U \X )^{- 1} ( \U \X )^{\top } \le \U \non 
\end{align}
as in (20) of Wang and Roy (2018). 
Therefore, 
\begin{align}
p( \bbe , \bSi | \u ) &\ge %
\{ c_1 / | \bSi |^{(n + \nu + d + 1 + p) / 2} \} \etr \{ - \bSi ^{- 1} ( \bTh ^{- 1} + 5 \B ^{\top } \A ^{- 1} \B + 2 \bBe ^{\top } \A ^{- 1} \bBe ) / 2 \} \non \\
&\quad \times \etr \{ - \bSi ^{- 1} (2 \bBe ^{\top } \X ^{\top } \U \X \bBe + 5 \Y ^{\top } \U \Y ) / 2 \} \non \\
&= c_1 {\etr \{ - \bSi ^{- 1} ( \bTh ^{- 1} + 5 \B ^{\top } \A ^{- 1} \B + 2 \bBe ^{\top } \A ^{- 1} \bBe ) / 2 \} \over | \bSi |^{(n + \nu + d + 1 + p) / 2}} \exp \Big( - {1 \over 2} \sum_{i = 1}^{n} ( \bXi )_{i, i} u_i \Big) \text{,} \non 
\end{align}
where $\bXi = 2 \X \bBe \bSi ^{- 1} \bBe ^{\top } \X ^{\top } + 5 \Y \bSi ^{- 1} \Y ^{\top } \ge \O ^{(n)}$. 
By the covariance inequality, for all $i = 1, \dots , n$, we have 
\begin{align}
&E^{\u | ( \bBe , \bSi )} [ \exp (- ( \bXi )_{i, i} u_i / 2) | \bBe ^{\rm{o}} , \bSi ^{\rm{o}} ] \non \\
&= \frac{ \displaystyle \int_{0}^{\infty } h( u_i ) {u_i}^{d / 2} \exp \{ - ( u_i / 2) ( \y _i - \bBe ^{\top } \x _i )^{\top } \bSi ^{- 1} ( \y _i - \bBe ^{\top } \x _i ) \} \exp (- ( \bXi )_{i, i} u_i / 2) d{u_i} }{ \displaystyle \int_{0}^{\infty } h( u_i ) {u_i}^{d / 2} \exp \{ - ( u_i / 2) ( \y _i - \bBe ^{\top } \x _i )^{\top } \bSi ^{- 1} ( \y _i - \bBe ^{\top } \x _i ) \} d{u_i} } \non \\
&\ge \int_{0}^{\infty } h( u_i ) {u_i}^{d / 2} \exp (- ( \bXi )_{i, i} u_i / 2) d{u_i} / \int_{0}^{\infty } h( u_i ) {u_i}^{d / 2} d{u_i} \text{,} \non 
\end{align}
where the denominator of the right-hand side is finite by assumption. 
Thus, 
\begin{align}
&E^{\u | ( \bBe , \bSi )} [ p( \bBe , \bSi | \u ) | \bBe ^{\rm{o}} , \bSi ^{\rm{o}} ] \non \\
&\ge \{ c_1 / | \bSi |^{(n + \nu + d + 1 + p) / 2} \} \etr \{ - \bSi ^{- 1} ( \bTh ^{- 1} + 5 \B ^{\top } \A ^{- 1} \B + 2 \bBe ^{\top } \A ^{- 1} \bBe ) / 2 \} \non \\
&\quad \times \prod_{i = 1}^{n} \Big\{ \int_{0}^{\infty } h( u_i ) {u_i}^{d / 2} \exp (- (2 \X \bBe \bSi ^{- 1} \bBe ^{\top } \X ^{\top } + 5 \Y \bSi ^{- 1} \Y ^{\top } )_{i, i} u_i / 2) d{u_i} / \int_{0}^{\infty } h( u_i ) {u_i}^{d / 2} d{u_i} \Big\} \text{.} \non 
\end{align}
Since the right hand side is an integrable function of $( \bBe , \bSi )$, the result follows. 
\hfill$\Box$

\section{Geometric Ergodicity When %
the Prior is Improper}
\label{sec:improper} 
In the last section, we saw that the Markov chain is geometrically ergodic even if the mixing density is heavy-tailed near the origin, in the case where we use a conditionally conjugate proper prior. 
A natural question then is whether the assumptions on the left tail of the mixing density which are assumed for improper cases in the literature are indeed necessary; see the discussion on page 517 of Hobert et al. (2018). 
In this section, we show that geometric ergodicity holds under significantly milder conditions in the simple case treated by Roy and Hobert (2010). 
Although we can always use a proper prior in practice, the following result of this section is relevant because improper priors can be interpreted as limiting cases of proper priors.

\begin{prp}
\label{prp:improper} 
Consider the data augmentation algorithm of Section 4 of Roy and Hobert (2010), which is used for the model with the multivariate Student's $t$ error distribution with $\nu > 0$ degrees of freedom and with the improper prior $\pi ( \bBe , \bSi ) \propto | \bSi |^{- (d + 1) / 2}$ for the parameters. 
Assume that $n \ge p + d$ so that the posterior is proper (by Proposition 1 of Roy and Hobert (2010)). 
Suppose that $( \X , \Y ) \in \mathbb{R} ^{n \times (p + d)}$ is of full rank and that $\nu + d > 2$. 
Suppose that there exists $\ep > 0$ such that 
\begin{align}
&{n - p \over \nu + d - 2} < 1 / \tr \Big\{ \Big( \sum_{i = 1}^{n} {\w _i {\w _i}^{\top } \over \ep + \| \w _i \| ^2}  \Big) ^{- 1}\Big\} \text{,} \label{eq:condition_improper} 
\end{align}
where $\w _i = ( {\x _i}^{\top } , {\y _i}^{\top } )^{\top } \in \mathbb{R} ^{p + d}$ for $i = 1, \dots , n$. 
Then the Markov chain based on their data augmentation algorithm is geometrically ergodic. 
\end{prp}

Theorem 1 of Roy and Hobert (2010) assumes that $n < \nu + p - 2$, which is restrictive when $n$ is larger. 
In contrast, (\ref{eq:condition_improper}) is expected to be less restrictive, for if $n$ is large, so will the right-hand side tend to be. 

In order to prove Proposition \ref{prp:improper}, we consider the following class of enegy functions: 
\begin{align}
V_{\C } ( \bBe , \bSi ) &= \tr \{ ( \Y - \X \bBe ) \bSi ^{- 1} ( \Y - \X \bBe )^{\top } \C \} \text{,} \quad \bBe \in \mathbb{R} ^{p \times d} \text{,} \quad \bSi > \O ^{(d)} \text{,} \non 
\end{align}
with $\C > \O ^{(n)}$. 
The energy function of Lemma 1 of Roy and Hobert (2010) is obtained by setting $\C = \I ^{(n)}$. 
We set $\C = \C _0$, where 
\begin{align}
\C _0 &= (( \X , \Y , \Z )^{- 1} )^{\top } ( \X , \Y , \Z )^{- 1} \non 
\end{align}
for some fixed $\Z \in \mathbb{R} ^{n \times (n - p - d)}$ such that $( \X , \Y , \Z ) \in \mathbb{R} ^{n \times n}$ is nonsingular. 
Then we have the simple quadratic form 
\begin{align}
V_{\C _0} ( \bBe , \bSi ) &= \tr ( \bSi ^{- 1} ) + \tr ( \bBe \bSi ^{- 1} \bBe ^{\top } ) \label{eq:energy_quadratic} 
\end{align}
since $\X ^{\top } \C _0 \X = \I ^{(p)}$, $\Y ^{\top } \C _0 \Y = \I ^{(d)}$, and $\X ^{\top } \C _0 \Y = \O ^{(p, d)}$. 

\bigskip

\noindent
{\bf Proof of Proposition \ref{prp:improper}.} \ \ Since $\C _0 \ge \I ^{(n)} / \tr ( {\C _0}^{- 1} )$, we have $V_{\I ^{(n)}} ( \bBe , \bSi ) \le \{ \tr ( {\C _0}^{- 1} ) \} V_{\C _0} ( \bBe , \bSi )$ and the minorization condition associated with (\ref{eq:energy_quadratic}) is veryfied as in the proof of Lemma 2 of Roy and Hobert (2010). 

We suppress the dependence on $\Y $ and use the notation $u_i$, $\U $, and $\bGa $ instead of $q_i$, $\Q ^{- 1}$, and $\bmu $ for $i = 1, \dots , n$. 
By results in the proof of Lemma 1 of Roy and Hobert (2010), 
\begin{align}
E[ V_{\C _0} ( \bBe , \bSi ) | \u ] &= \tr ( E[ \bSi ^{- 1} | \u ]) + \tr \{ d ( \X ^{\top } \U \X )^{- 1} + \bGa E[ \bSi ^{- 1} | \u ] \bGa ^{\top } \} \non 
\end{align}
and 
\begin{align}
E[ \bSi ^{- 1} | \u ] &= (n - p) \{ \Y ^{\top } \U \Y - \bGa ^{\top } ( \X ^{\top } \U \X ) \bGa \} ^{- 1} \non \\
&= (n - p) \Big\{ \sum_{i = 1}^{n} ( \y _i - \bGa ^{\top } \x _i ) u_i ( \y _i - \bGa ^{\top } \x _i )^{\top } \Big\} ^{- 1} \non \\
&= (n - p) \{ ( \Y - \X \bGa )^{\top } \U ( \Y - \X \bGa ) \} ^{- 1} \text{.} \non 
\end{align}
Therefore, 
\begin{align}
E[ V_{\C _0} ( \bBe , \bSi ) | \u ] &= d \tr \{ ( \X ^{\top } \U \X )^{- 1} \} + (n - p) \non \\
&\quad \times ( \tr [ \{ ( \Y - \X \bGa )^{\top } \U ( \Y - \X \bGa ) \} ^{- 1} ] + \tr [ \bGa \{ ( \Y - \X \bGa )^{\top } \U ( \Y - \X \bGa ) \} ^{- 1} \bGa ^{\top } ]) \text{.} \non 
\end{align}

Let $\Xbt = \U ^{1 / 2} \X $ and $\Ybt = \U ^{1 / 2} \Y $ and write 
\begin{align}
&\P _{\M } = \M ( \M ^{\top } \M )^{+} \M ^{\top } \quad \text{and} \quad \Q _{\M } = \I ^{(n)} - \P _{\M } \text{,} \non 
\end{align}
which are symmetric and idempotent, for any matrix $\M $ whose number of rows is $n$. 
Then 
\begin{align}
&\{ ( \Y - \X \bGa )^{\top } \U ( \Y - \X \bGa ) \} ^{- 1} = \{ ( \Q _{\Xbt } \Ybt )^{\top } ( \Q _{\Xbt } \Ybt ) \} ^{- 1} = ( \Ybt ^{\top } \Q _{\Xbt } \Ybt )^{- 1} \text{.} \non 
\end{align}
Note that 
\begin{align}
( \Ybt ^{\top } \Q _{\Xbt } \Ybt )^{- 1} &= \{ \Ybt ^{\top } \Ybt - \Ybt ^{\top } \Xbt ( \Xbt ^{\top } \Xbt )^{- 1} \Xbt ^{\top } \Ybt \} ^{- 1} \non \\
&= ( \Ybt ^{\top } \Ybt )^{- 1} - ( \Ybt ^{\top } \Ybt )^{- 1} \Ybt ^{\top } \Xbt \{ - \Xbt ^{\top } \Xbt + \Xbt ^{\top } \Ybt ( \Ybt ^{\top } \Ybt )^{- 1} \Ybt ^{\top } \Xbt \} ^{- 1} \Xbt ^{\top } \Ybt ( \Ybt ^{\top } \Ybt )^{- 1} \non \\
&= ( \Ybt ^{\top } \Ybt )^{- 1} + ( \Ybt ^{\top } \Ybt )^{- 1} \Ybt ^{\top } \Xbt ( \Xbt ^{\top } \Q _{\Ybt } \Xbt )^{- 1} \Xbt ^{\top } \Ybt ( \Ybt ^{\top } \Ybt )^{- 1} \text{.} \non 
\end{align}
Then 
\begin{align}
&\bGa \{ ( \Y - \X \bGa )^{\top } \U ( \Y - \X \bGa ) \} ^{- 1} \bGa ^{\top } \non \\
&= ( \Xbt ^{\top } \Xbt )^{- 1} \Xbt ^{\top } \P _{\Ybt } \Xbt ( \Xbt ^{\top } \Xbt )^{- 1} + ( \Xbt ^{\top } \Xbt )^{- 1} \Xbt ^{\top } \P _{\Ybt } \Xbt ( \Xbt ^{\top } \Q _{\Ybt } \Xbt )^{- 1} \Xbt ^{\top } \P _{\Ybt } \Xbt ( \Xbt ^{\top } \Xbt )^{- 1} \non \\
&= ( \Xbt ^{\top } \Xbt )^{- 1} \Xbt ^{\top } \P _{\Ybt } \Xbt ( \Xbt ^{\top } \Xbt )^{- 1} - ( \Xbt ^{\top } \Xbt )^{- 1} \Xbt ^{\top } ( \Q _{\Ybt } - \I ^{(n)} ) \Xbt ( \Xbt ^{\top } \Q _{\Ybt } \Xbt )^{- 1} \Xbt ^{\top } \P _{\Ybt } \Xbt ( \Xbt ^{\top } \Xbt )^{- 1} \non \\
&= ( \Xbt ^{\top } \Xbt )^{- 1} \Xbt ^{\top } \P _{\Ybt } \Xbt ( \Xbt ^{\top } \Xbt )^{- 1} - ( \Xbt ^{\top } \Xbt )^{- 1} \Xbt ^{\top } \P _{\Ybt } \Xbt ( \Xbt ^{\top } \Xbt )^{- 1} + ( \Xbt ^{\top } \Q _{\Ybt } \Xbt )^{- 1} \Xbt ^{\top } \P _{\Ybt } \Xbt ( \Xbt ^{\top } \Xbt )^{- 1} \non \\
&= ( \Xbt ^{\top } \Q _{\Ybt } \Xbt )^{- 1} \Xbt ^{\top } \P _{\Ybt } \Xbt ( \Xbt ^{\top } \Xbt )^{- 1} = ( \Xbt ^{\top } \Q _{\Ybt } \Xbt )^{- 1} - ( \Xbt ^{\top } \Xbt )^{- 1} \text{.} \non 
\end{align}
Therefore, 
\begin{align}
&[E[ V_{\C _0} ( \bBe , \bSi ) | \u ] - d \tr \{ ( \X ^{\top } \U \X )^{- 1} \} / (n - p) \non \\
&= \tr \{ ( \Ybt ^{\top } \Q _{\Xbt } \Ybt )^{- 1} \} + \tr \{ ( \Xbt ^{\top } \Q _{\Ybt } \Xbt )^{- 1} - ( \Xbt ^{\top } \Xbt )^{- 1} \} \non 
\end{align}
or 
\begin{align}
&E[ V_{\C _0} ( \bBe , \bSi ) | \u ] = (n - p) [ \tr \{ ( \Ybt ^{\top } \Q _{\Xbt } \Ybt )^{- 1} \} + \tr \{ ( \Xbt ^{\top } \Q _{\Ybt } \Xbt )^{- 1} \} ] - (n - p - d) \tr \{ ( \Xbt ^{\top } \Xbt )^{- 1} \} \text{.} \non 
\end{align}
Thus, 
\begin{align}
E[ V_{\C _0} ( \bBe , \bSi ) | \u ] &\le (n - p) [ \tr \{ ( \Ybt ^{\top } \Q _{\Xbt } \Ybt )^{- 1} \} + \tr \{ ( \Xbt ^{\top } \Q _{\Ybt } \Xbt )^{- 1} \} ] \non \\
&= (n - p) \tr \Big( \begin{pmatrix} \Xbt ^{\top } \Xbt & \Xbt ^{\top } \Ybt \\ \Ybt ^{\top } \Xbt & \Ybt ^{\top } \Ybt \end{pmatrix} ^{- 1} \Big) \non \\
&= (n - p) \tr \Big\{ \Big( \begin{pmatrix} \X ^{\top }\\ \Y ^{\top } \end{pmatrix} \U \begin{pmatrix} \X & \Y \end{pmatrix} \Big) ^{- 1} \Big\} \text{,} \label{pimproperp1} 
\end{align}
where the first equality follows, for example, from Theorem A.3.3 of Anderson (2003). 

Now, let $\w _i = ( {\x _i}^{\top } , {\y _i}^{\top } )^{\top } \in \mathbb{R} ^{p + d}$ for $i = 1, \dots , n$. 
For each $i = 1, \dots , n$, 
\begin{align}
\tr \Big\{ \Big( \begin{pmatrix} \X ^{\top } \\ \Y ^{\top } \end{pmatrix} \U \begin{pmatrix} \X & \Y \end{pmatrix} \Big) ^{- 1} \Big\} &= \tr \{ ( \M _i + \w _i u_i {\w _i}^{\top } ) ^{- 1} \} \non \\
&= \tr \{ {\M _i}^{- 1} - {\M _i}^{- 1} \w _i (1 / u_i + {\w _i}^{\top } {\M _i}^{- 1} \w _i ) ^{- 1} {\w _i}^{\top } {\M _i}^{- 1} \} \non \\
&= \tr ( {\M _i}^{- 1} ) - {{\w _i}^{\top } {\M _i}^{- 1} {\M _i}^{- 1} \w _i \over 1 / u_i + {\w _i}^{\top } {\M _i}^{- 1} \w _i} \non 
\end{align}
is a concave function of $1 / u_i$, where 
\begin{align}
\M _i &= \sum_{\substack{j = 1 \\ j \neq i}}^{n} \w _j u_j {\w _j}^{\top } \text{.} \non 
\end{align}
Therefore, by Jensen's inequality, 
\begin{align}
E \Big[ \tr \Big\{ \Big( \begin{pmatrix} \X ^{\top } \\ \Y ^{\top } \end{pmatrix} \U \begin{pmatrix} \X & \Y \end{pmatrix} \Big) ^{- 1} \Big\} \Big| \bBe , \bSi \Big] &\le \tr \Big[ \Big\{ \begin{pmatrix} \X ^{\top }\\ \Y ^{\top } \end{pmatrix} (E[ \U ^{- 1} | \bBe , \bSi ])^{- 1} \begin{pmatrix} \X & \Y \end{pmatrix} \Big\} ^{- 1} \Big] \text{.} \non 
\end{align}
Expressions for the conditional means of ${u_i}^{- 1}$, $i = 1, \dots , n$, are given in the proof of Lemma 1 of Roy and Hobert (2010) and we have 
\begin{align}
&E[ {u_i}^{- 1} | \bBe , \bSi ] = {\nu + {\w _i}^{\top } \bTh \w _i \over \nu + d - 2} \le {\nu + ( \tr \bTh ) \| \w _i \| ^2 \over \nu + d - 2} = {\nu + \| \w _i \| ^2 V_{\C _0} ( \bBe , \bSi ) \over \nu + d - 2} \non 
\end{align}
for all $i = 1, \dots , n$, where 
\begin{align}
\bTh &= \begin{pmatrix} \bBe \bSi ^{- 1} \bBe ^{\top } & - \bBe \bSi ^{- 1} \\ - \bSi ^{- 1} \bBe ^{\top } & \bSi ^{- 1} \end{pmatrix} \text{.} \non 
\end{align}
Thus, 
\begin{align}
E \Big[ \tr \Big\{ \Big( \begin{pmatrix} \X ^{\top } \\ \Y ^{\top } \end{pmatrix} \U \begin{pmatrix} \X & \Y \end{pmatrix} \Big) ^{- 1} \Big\} \Big| \bBe , \bSi \Big] &\le \tr \Big[ \Big\{ \sum_{i = 1}^{n} \w _i {\nu + d - 2 \over \nu + \| \w _i \| ^2 V_{\C _0} ( \bBe , \bSi )} {\w _i}^{\top } \Big\} ^{- 1} \Big] \non \\
&= {V_{\C _0} ( \bBe , \bSi ) \over \nu + d - 2} \tr \Big[ \Big\{ \sum_{i = 1}^{n} {\w _i {\w _i}^{\top } \over \nu / V_{\C _0} ( \bBe , \bSi ) + \| \w _i \| ^2} \Big\} ^{- 1} \Big] \text{.} \label{pimproperp2} 
\end{align}

Fix $\bBe ^{\rm{o}} \in \mathbb{R} ^{p \times d}$ and $\bSi ^{\rm{o}} > \O ^{(d)}$. 
By (\ref{pimproperp1}) and (\ref{pimproperp2}), 
\begin{align}
&E[ E[ V_{\C _0} ( \bBe , \bSi ) | \u ] | ( \bBe , \bSi ) = ( \bBe ^{\rm{o}} , \bSi ^{\rm{o}} ) ] \non \\
&\le (n - p) {V_{\C _0} ( \bBe ^{\rm{o}} , \bSi ^{\rm{o}} ) \over \nu + d - 2} \tr \Big[ \Big\{ \sum_{i = 1}^{n} {\w _i {\w _i}^{\top } \over \nu / V_{\C _0} ( \bBe ^{\rm{o}} , \bSi ^{\rm{o}} ) + \| \w _i \| ^2} \Big\} ^{- 1} \Big] \text{.} \non 
\end{align}
For any $M_1 > 0$, 
\begin{align}
&E[ E[ V_{\C _0} ( \bBe , \bSi ) | \u ] | ( \bBe , \bSi ) = ( \bBe ^{\rm{o}} , \bSi ^{\rm{o}} ) ] \le (n - p) {M_1 \over \nu + d - 2} \tr \Big\{ \Big( \sum_{i = 1}^{n} {\w _i {\w _i}^{\top } \over \nu / M_1 + \| \w _i \| ^2} \Big) ^{- 1} \Big\} \non 
\end{align}
if $V_{\C _0} ( \bBe ^{\rm{o}} , \bSi ^{\rm{o}} ) \le M_1$ and 
\begin{align}
&E[ E[ V_{\C _0} ( \bBe , \bSi ) | \u ] | ( \bBe , \bSi ) = ( \bBe ^{\rm{o}} , \bSi ^{\rm{o}} ) ] \le (n - p) {V_{\C _0} ( \bBe ^{\rm{o}} , \bSi ^{\rm{o}} ) \over \nu + d - 2} \tr \Big\{ \Big( \sum_{i = 1}^{n} {\w _i {\w _i}^{\top } \over \nu / M_1 + \| \w _i \| ^2} \Big) ^{- 1} \Big\} \non 
\end{align}
if $V_{\C _0} ( \bBe ^{\rm{o}} , \bSi ^{\rm{o}} ) > M_1$. 
Since $\nu / M_1 < \ep $ for some $M_1 > 0$, the drift condition holds by (\ref{eq:condition_improper}). 
This completes the proof. 
\hfill$\Box$

\section*{Acknowledgments}
Research of the author was supported in part by JSPS KAKENHI Grant Number JP22K20132 from Japan Society for the Promotion of Science.

\end{document}